\documentclass[12pt]{amsart}

\usepackage{amssymb}
\usepackage{bbm}
\usepackage{enumerate, xspace}

\let\mathds\mathbf

\newcommand{\one}{{\mathds{1}}}

\renewcommand{\d}{\;d\:}

\newcommand{\Xjtruncatedeps}{X_j \one_{\{0< |X_j|< b_n \epsilon\}}}
\renewcommand{\Xjtruncatedeps}{X_j \one_{\{|X_j|< b_n \epsilon\}}}

\newcommand{\Xjtruncatedone}{X_j \one_{\{0< |X_j|< b_n\}}}
\renewcommand{\Xjtruncatedone}{X_j \one_{\{|X_j|< b_n\}}}

\let\epsilon\varepsilon

\newcommand\ds{\displaystyle}

\newcommand\dto{\, {\stackrel{d}{\to} \,}}

\newcommand{\Lebesgue}{\operatorname{Leb}}

\newcommand{\BV}{\operatorname{BV}}

\newcommand{\bvnorm}[1]{\|#1\|_{\operatorname{BV}}}

\newcommand{\sign}[1]{{\text{sign}}(#1)}

\numberwithin{equation}{section}

\newtheorem{thm}{Theorem}[section]
\newtheorem{lem}[thm]{Lemma}
\newtheorem{cor}[thm]{Corollary}

\newtheorem{prop}[thm]{Proposition}

\newtheorem{defn}[thm]{Definition}

\newtheorem{rmk}[thm]{Remark}

\newcommand\cB{{\mathcal B}}

\newcommand\cP{{\mathcal P}}

\newcommand\bE{{\mathbb E}}

\newcommand\bN{{\mathbb N}}

\newcommand\bR{{\mathbb R}}

\newcommand{\R}{\mathbb{R}}

\begin{document}

\date{\today.}

\title[Birkhoff sum convergence of Fr\'echet observables to stable
laws]{Birkhoff sum convergence of Fr\'echet observables to stable laws for
  Gibbs-Markov systems and applications.}

\author[A.Chen]{An Chen}
\address{An Chen\\ Department of Mathematics\\
  University of Houston\\
  Houston\\
  TX 77204\\
  USA} \email{anchenwhu@gmail.com}

\author[M. Nicol]{Matthew Nicol}
\address{Matthew Nicol\\ Department of Mathematics\\
University of Houston\\
Houston\\
TX 77204\\
USA} \email{nicol@math.uh.edu}
\urladdr{http://www.math.uh.edu/~nicol/}

\author[A. T\"or\"ok]{Andrew T\"or\"ok}
\address{Andrew T\"or\"ok\\ Department of Mathematics\\
  University of Houston\\
  Houston\\
  TX 77204\\
  USA and 
{Institute of Mathematics of the Romanian Academy, Bucharest, Romania.}}
\email{torok@math.uh.edu}
\urladdr{http://www.math.uh.edu/~torok/}

\thanks{MN was supported in
  part by NSF Grant DMS 2009923 and would like to Jorge Freitas and Roland Zweim\"uller for helpful discussions
  and  the Erwin Schr\"{o}dinger Institute for support through the 2024 Workshop on Extremes and Rare Events. AT was supported in part
  by NSF Grant DMS 1816315. }



\begin{abstract}

  We use a Poisson point process approach to prove distributional convergence to a stable law for
  non square-integrable observables $\phi: [0,1]\to \bR$, mostly of the
  form $\phi (x) = d(x,x_0)^{-\frac{1}{\alpha}}$,$0<\alpha\le 2$, on Gibbs-Markov maps.  A key result is to verify a standard mixing condition, which ensures that large values of the observable dominate the time-series, in the range $1<\alpha \le 2$. Stable limit laws for observables on dynamical systems have been
established in two  settings: ``good observables''
(typically H\"older) on slowly mixing non-uniformly hyperbolic systems and
``bad'' observables (unbounded with fat tails) on fast mixing dynamical
systems. As an application we investigate the interplay between these two effects in a class of intermittent-type maps.

  \end{abstract}

\maketitle

\tableofcontents

\section{Introduction}\label{sec:intro}

In this paper we consider distributional convergence to  stable laws for
non square-integrable observables $\phi: [0,1]\to \bR$ of form
$\phi (x) = d(x,x_0)^{-\frac{1}{\alpha}}$,$0<\alpha\le 2$, on Gibbs-Markov
maps of the unit interval $[0,1]$ ($x_0 \in [0,1]$). Our results imply distributional convergence, in some parameter regimes, to  stable laws for
non square-integrable observables on certain systems modeled by first return time Young Towers in
which the base map is Gibbs-Markov, in particular intermittent-type maps of the unit interval.

Most of our results consider distance-like observables
$\phi (x) = d(x,x_0)^{-\frac{1}{\alpha}}$, where $\alpha\in (0,2)$ and
$x_0\in [0,1]$. But  our result on mixing conditions, Theorem~\ref{thm:mix_NC},  extends to observables 
$\phi$ which are  regularly varying with stable  index $\alpha$ and for which, for sufficiently large $t$, $\|\phi\one_{|\phi|<t}\|_{BV} \le Kt$ for some constant $K$.

Stable limit laws for observables on dynamical systems have been
established in two somewhat distinct settings: ``good observables''
(typically H\"older) on slowly mixing non-uniformly hyperbolic systems and
``bad'' observables (unbounded with fat tails) on fast mixing dynamical
systems.

For results on the first type we refer to the influential
papers~\cite{Gouezel_Intermittent,Gouezel_Skew}
and~\cite{Melbourne_Zweimuller}. In the setting of ``good observables''
(typically H\"older) on slowly mixing non-uniformly hyperbolic systems the
technique of inducing on a subset of phase space and constructing a Young
Tower has been used with some success. ``Good'' observables lift to
well-behaved observables lying in a suitable Banach space on the Young
Tower. This is not the case in general  with unbounded observables with fat tails,
though in~\cite{Gouezel_Intermittent} the induction technique permits analysis of an
observable which is unbounded at the fixed point $x=0$ in a family of
intermittent maps. As $x=0$ is not in the Young Tower the observable lifts to a function on the 
Tower which is bounded on each  column of the Tower and  with sufficient regularity for spectral techniques to apply.

For general  results on distributional and functional stable laws for non-square
integrable observables using a Poisson point process approach we refer to the papers of  Marta
Tyran-Kaminska~\cite{TK,TK-dynamical}.  Tyran-Kaminska considers convergence of Birkhoff sums to stable laws
and corresponding functional convergence in the $J_1$ topology to L\'evy processes. She uses a point process approach but
her work explicitly excludes clustering behavior, and so is not applicable to observables $\phi (x) = d(x,x_0)^{-\frac{1}{\alpha}}$
maximized at a periodic point $x_0$ (for which clustering of extremes is expected).
In the setting of Gibbs-Markov maps
Tyran-Kaminska  shows, among other results, that functions which are measurable with
respect to the Gibbs-Markov partition and in the domain of attraction of a
stable law with index $\alpha$ converge (under the appropriate scaling) in
the $J_1$ topology to a L\'evy process of index $\alpha$~\cite[Theorem 3.3,
Corollaries 4.1 and 4.2]{TK-dynamical}. Her result is not applicable in our
setting as $\phi (x) = d(x,x_0)^{-\frac{1}{\alpha}}$ is not measurable with
respect to the Gibbs-Markov partition in the settings we consider. It is
interesting to note though that in the case of a slowly-mixing intermittent
map with an indifferent fixed point at $x=0$ and a H\"older observable
$\phi$, $\phi(0)\not =0$, the constant $\phi (0)$ may be induced as a
measurable function on the Gibbs-Markov base of the usual first return
tower representation. This approach is used   by Melbourne and Zweimuller~\cite{Melbourne_Zweimuller}
to prove convergence to stable laws for H\"older functions on slowly-mixing systems modeled by 
a Young Tower.

Marta Tyran-Kaminska's work is based on a Poisson point process approach described by Durrett
and Resnick~\cite{Resnick1,Resnick2,Durrett_Resnick}. This paper follows a
similar approach to the scheme laid out by Tyran-Kaminska in that we
require convergence of a counting process to a  Poisson
process and a form of decay of correlations estimate for  a truncation of the  observable that ensures the
Birkhoff sum of small values of $\phi$ do not contribute too much and it is
the large values that dominate. We stress that we do not prove functional convergence in
the $J_1$ topology but rather distributional convergence. In fact, as
Tyran-Kaminska shows in~\cite[Theorem 1.1, Example 1.1]{TK-dynamical} in
situations where the counting process exhibits clustering convergence  in
the $J_1$ topology does not hold. Recent work has shown that in some  settings where $J_1$ convergence 
does not hold  that
convergence is possible in the weaker $M_1$
topology~\cite{Melbourne_Zweimuller}  and in the  $F^{'}$ topology~\cite{FFTV23}. We refer to
these papers for helpful discussions of the relevant topologies and related
results.

In the cases where we obtain  distributional rather than functional convergence, we
need only validate  the weaker conditions of Davis and Tsing~\cite[Theorem 3.1]{DH95} rather than the stronger condition of ~\cite[Theorem 1.1 Condition (1.5)]{TK-dynamical}. This allows us to consider a different class of  observables than in~\cite{TK-dynamical}.

In the setting of Gibbs-Markov maps (or more generally Rychlik maps)
Freitas, Freitas and Magalhaes~\cite{FFM20} have proved that observables of the
type $d(x,x_0)=d(x,x_0)^{-\frac{1}{\alpha}}$, $x_0\not =0$, have counting
processes that converge to a simple Poisson point process if $x_0$ is not
periodic and a  ``clustered' point  process if $x_0$ is periodic. Furthermore 
if $0<\alpha<1$ then~\cite{FFT20} have shown functional convergence of  the rescaled
time-series for this observable in the $F^{'}$ topology, which implies convergence of the 
scaled Birkhoff sum to a stable law. One contribution of this
paper is  Theorem~\ref{thm:mix_NC} which verifies a mixing condition in the case $1<\alpha<2$ and extends
the stable law convergence to the parameter range $1<\alpha <2$.

One question that arose in our investigation (that was not satisfactorily resolved) can be stated simply. Suppose ($T_{\gamma},[0,1],\mu_{\gamma})$ is a LSV~\cite{LSV99} map of the 
unit interval (see Section 6) and $\mu_{\gamma}$ is the Lebesgue equivalent invariant measure for $T_{\gamma}$. 
Suppose $\phi$ has support in $[1/2,1]$, $\int \phi\; d {\mu_{\gamma}}=0$,  and locally, near $x_0\in [1/2,1]$ is of form $d(x,x_0)^{-\frac{1}{\alpha}}$ (elsewhere H\"older). We are  able to show that the Birkhoff sum of the induced map on $[1/2,1]$ converges in distribution to  a stable law with index $\alpha $. In  certain parameter regions for $1\le \alpha \le2$ and $0<\gamma <1$, namely $\frac{1}{\gamma}\le \alpha \le 1+\frac{1}{\gamma^2}-\frac{1}{\gamma}$, we are able to show that the stable law with index $\alpha$ lifts  from that of the induced observable to  give a stable law for the 
original observable $\phi$.  Does a  stable law  of index $\alpha$ lift for all $\alpha<\frac{1}{\gamma}$  for all parameter ranges of  $1 \le \alpha \le 2$ and $0<\gamma <1$ if $\int \phi \;d{\mu_{\gamma}}=0$ and $\phi$ has support in $[1/2,1]$?


Our main results are given in the section~\ref{Main_Results}.  We first give some
background.

 \section{Probabilistic tools}\label{sec:tools}

In this section, we review some topics from Probability Theory.

\subsection{Regularly varying functions and domains of
  attraction}\label{ssec:regu}

We refer to
Feller~\cite{Fel71} or Bingham, Goldie and
Teugels~\cite{Bingham-Goldie-Teugels-1987}  for the relations between domains of attraction of
stable laws and regularly varying functions.  For $\phi$ regularly varying  we define scaling constants $b_n$ (related to the index) and 
$c_n$ (centering) by
\begin{defn}\label{def:scaling-constants}
  Given a regularly varying function $\phi$  of index $\alpha\in (0,2)$ on a probability space
  $(X, \mu)$, introduce:

  \begin{enumerate}[--]
  \item a scaling sequence $(b_n)_{n \ge 1}$ by
    \begin{equation}\label{eq:tail1}
      \lim_{n \to \infty} n \mu( |\phi|  > b_n ) = 1.
    \end{equation}

  \item a centering sequence $(c_n)_{n\ge 1}$ by
    \begin{equation}\label{eq:centering}
      c_n =
      \begin{cases}
        0 &\text{if $\alpha \in (0,1)$}\cr n 
        E [\phi ] & \text{if $\alpha \in(1,2)$}
      \end{cases}.
    \end{equation}
  \end{enumerate}
\end{defn}

The constants $p,q$ are defined by  
\[
p=\lim_{t\to \infty} \frac{\mu(\phi >t)}{\mu(|\phi|>t)}
\]
and $q=1-p$.

Note that if  $\phi= d(x,x_0)^{ - \frac{1}{\alpha}}$ is an observable on the unit interval $[0,1]$
equipped with a Lebesgue equivalent measure and $x_0\in [0,1]$  then 
$b_n\sim n^{1/\alpha}$, where $\sim$ means there exists $C_1$, $C_2>0$ with $C_1 n^{1/\alpha} \le b_n \le C_2 n^{1/\alpha}$. Note also
that  $p=1$ as $\phi > 0$. As we did above, we will sometimes write $E[\phi]$ for the
expectation of an observable when the measure is clear from
context.

\begin{rmk}\label{rem:centering-constants}
  When $\alpha \in (0,1)$ then $\phi$ is not integrable and one can choose
  the centering sequence $(c_n)$ to be identically $0$. When $\alpha = 1$,
  it might happen that $\phi$ is not integrable, and it is then necessary to center. We don't consider the case $\alpha=1$. In the literature if centering is needed it is often 
  specified as $c_n=n\mathbb{E}(\phi\one_{|\phi|\le b_n})$ but we have opted for a simpler scaling.
   By \cite[Remark 3.1]{DH95}, for $1<\alpha \le 2$ if  $\phi$ is a regularly varying function of index $\alpha$
  then  $nE(\phi)$ may be   used in scaling rather than the truncation in~\eqref{eq:stable-alpha-ge1}. The same limiting distribution $S$ is obtained  though shifted by the constant $(p-q)\frac{\alpha}{\alpha-1}$.
  More precisely
      \begin{equation*} 
      \frac{1}{b_n} \left(\sum_{j=1}^{n} [\phi \circ T^j - \mu (\phi)]
      \right) \to_{d} S-(p-q)\frac{\alpha}{\alpha-1}
    \end{equation*}
    where $q=1-p$.
    This is a consequence of 
    \[
    \frac{n}{b_n} [E(\phi)-E(\phi\one_{\{|\phi
          |<b_n\}})]=\frac{n}{b_n} E[\phi\one_{(b_n,\infty)}(|\phi|)]\to (p-q)\frac{\alpha}{\alpha-1}
          \]
          using Karamata, so by convergence of types
          \[
           \frac{1}{b_n} \left(\sum_{j=1}^{n} \phi \circ T^j - c_n
      \right) \to_{d} S-(p-q)\frac{\alpha}{\alpha-1}
      \]

\end{rmk}

We will use the following asymptotics for truncated moments, which can be
deduced from Karamata's results concerning the tail behavior of regularly
varying functions.  Define $p$ by $ \lim_{x \to \infty}
    \frac{\nu(\phi>x)}{\nu(|\phi|>x)} = p$.

\begin{prop}[Karamata]\label{prop:karamata}

  Let $\phi$ be regularly varying with index $\alpha \in (0,2)$. Denote
  $\beta := 2p -1$ and, for $\epsilon > 0$,
  \begin{equation}\label{eq:c_alpha}
    c_\alpha(\epsilon) :=
    \begin{cases}
      0 &\text{if $\alpha \in (0,1)$}\cr
      - \beta \log \epsilon & \text{if $\alpha = 1$} \cr
      \epsilon^{1 - \alpha} \beta \alpha / (\alpha - 1) & \text{if
        $\alpha \in(1,2)$}
    \end{cases}
  \end{equation}
  The following hold for all $\epsilon>0$:\begin{enumerate}[(a)]

\item
  $\displaystyle \lim_{n \to \infty} n \mu ( |\phi| > \epsilon b_n ) =
  \epsilon^{-\alpha}$ (from the definition of $b_n$ and the regular
  variation of $\phi$)

\item If $k > \alpha$ then
  $\displaystyle E(|\phi|^k \one_{|\phi|\le u}) \sim
  \frac{\alpha}{k-\alpha}u^k \mu (|\phi|>u)$ as $u\to\infty$; in particular:


\item if $\alpha\in(0,2)$ then $E(| \phi|^2 \one_{\left\{ | \phi | \le \epsilon b_n \right\}}) \sim \frac{\alpha}{2 - \alpha} (\epsilon b_n)^2 \mu(| \phi | > \epsilon b_n),$
\item if $\alpha \in (0,1)$, $$E(| \phi| \one_{\left\{ | \phi | \le \epsilon b_n \right\}}) \sim \frac{\alpha}{1 - \alpha} \epsilon b_n \mu(| \phi | > \epsilon b_n),$$
\item if $\alpha \in (1,2)$, $$\lim_{n \to \infty} \frac{n}{b_n} E (\phi \one_{\{ | \phi | > \epsilon b_n \}}) = c_\alpha(\epsilon),$$
\item if $\alpha = 1$, $$\lim_{n \to \infty} \frac{n}{b_n} E (\phi \one_{\{  \epsilon b_n < |\phi | \le  b_n \}}) = c_\alpha(\epsilon),$$
\item if $\alpha = 1$, $$\frac{n}{b_n} E (|\phi| \one_{ \{  | \phi |
    \le \epsilon b_n \}}) \sim \widetilde{L}(n),$$ for a slowly varying
  function $\widetilde{L}$,
\end{enumerate}
\end{prop}


\subsection{L\'evy $\alpha$-stable processes}\label{ssec:stable}

A more detailed discussion of L\'evy processes is given in~\cite{TK, TK-dynamical}.

$X(t)$ is a L\'evy stable process if $X(0)=0$, $X$ has stationary independent increments and $X(1)$ has an 
$\alpha$-stable distribution.  Recall that the  distribution $F$  of a random variable $X$ is called $\alpha$-stable if 
there are constants $\gamma_n$ such  that for each $n$, if $X_i$ are iid with distribution $F$ then
\[
\sum_{j=1}^n X_j +\gamma_n \sim n^{\frac{1}{\alpha}} X_1
\]

The L\'evy-Khintchine representation for the characteristic function of an
$\alpha$-stable random variable $X_{\alpha,\beta}$ with index
$\alpha \in (0,2)$ and parameter $\beta\in[-1,1]$ has the form:
\[
\bE[e^{itX}]=\mbox{exp}\left[ita_{\alpha}+\int (e^{itx}-1-itx1_{[-1,1]} (x))\Pi_{\alpha}(dx)\right]
\]
where 
\begin{itemize}
\item 

  $\ds
  a_{\alpha}= \left\{ \begin{array}{ll}   \beta\frac{\alpha}{1-\alpha} & \mbox{ $\alpha \not = 1$}\\
                    0 & \mbox{ $\alpha =1$}\end{array} \right.,$

\item  $\Pi_{\alpha}$ is a L\'evy measure given by
  \[
    d \Pi_{\alpha} =\alpha (p1_{(0,\infty)} (x)+(1-p)1_{(-\infty,0)} (x) )
    |x|^{-\alpha -1} dx
 \]

\item $ \ds p=\frac{\beta+1}{2}.$
\end{itemize}

Note that $p$ and $\beta$ may equally serve as parameters for $X_{\alpha,\beta}$. We will drop the $\beta$ from $X_{\alpha,\beta}$, as is common in the literature, for simplicity of notation and when it plays no essential role.

\section{Stable law convergence}

Let $T$ be a measure preserving  transformation of a 
probability space $(X, \mu,\mathcal{B} )$.

Given $\phi : X \to \R$ measurable, we define the scaled
Birkhoff sum by
\begin{equation}\label{eq.sequential}
  S_n:= \frac{1}{b_n} [\sum_{j=0}^{ n - 1} \phi \circ
  T^j -  c_n],
\end{equation}
for some real constants $b_n > 0$, $c_n$.

We say $S_n$ converges to a stable law of index $\alpha$ if  
\[
S_n \dto X_{\alpha}
\]
for some random variable $X_{\alpha}$ with an $\alpha$-stable distribution.

\section{Poisson point processes}

Suppose $\phi$ is an observable on a dynamical system $(T,X,\mu)$ with stable index $\alpha$ and scaling constants $b_n$ and $c_n$. 
Let $B\subset (0,\infty) \times \R\setminus\{0\}$.

Define the counting process 
\[
N_n=\# \{ (j,(\phi \circ T^{j-1}-c_n)/b_n \in B\}
\]
For each $x\in (X,\mu)$, $N_n (x) $ is an integer valued counting process on $ (0,\infty) \times \R\setminus\{0\}$.

In our setting of Gibbs-Markov maps, Freitas, Freitas and Magalhaes~\cite{FFM20} have proved convergence of the counting measure $N_n$ 
( for $(T,X,\mu)$ a
Gibbs-Markov map and $\phi(x)=d(x,x_0)^{-1/\alpha}$) to a Poisson process which has the general form of~\cite[Corollary 2.4]{DH95}, namely
\[
N=\sum_{i=0}^{\infty} \sum_{j=1}^{\infty} \delta_{P_iQ_{ij}}
\]
where $\sum_{j=1}^{\infty} \delta_{P_i}$ is a Poisson process with intensity measure $\Pi_{\alpha}$ and $\sum_{j=1}^{\infty}\delta_{Q_{ij}}$ are point processes taking values
in $[-1,1]$ distributed according to a measure $\nu$. All point processes are mutually independent. In a dynamical setting, to which this Poisson point process is 
well suited, the $Q_{ij}$'s represent the ``clustering'' around an exceedance $P_i$ (which is chosen to be the largest value in the cluster).

\section{Gibbs-Markov Maps.}\label{GM}

We consider the following class of ergodic maps of $X=[0,1]$. Let $m$ denote Lebesgue measure and let $\mu$ be a Lebesgue equivalent measure
with density bounded above and away from zero below. Let $\cP$ 
be a countable partition of $[0,1]$ (mod $m$) into open intervals. 

 We suppose that
all partition elements $A_i \in \cP$ have $m(A_i) > 0$. A
$\mu$ measure-preserving transformation $T$ on $X$ is a Gibbs-Markov map if

(1) $\mathcal{B}$ is the smallest $\sigma$-algebra which contains $\bigvee_{n\geq 0} T^{-n}\cP$ which is complete with respect to $m$;

(2) Markov property: $\forall A_i \in \cP,$ $TA_i$ consists of a union of
partition elements and there exists $C>0$ such that $m(T A_i) >C$ for all
$i$. If $T:A_i \to X$ is onto $X$ mod $m$ for all $i$, we say that $T$ has
``full branches''.

(3) Local invertibility: $\forall A_i\in \cP$, $T:A_i \to TA_i$ is invertible. 


(4) Expansitivity:  There exists $\Lambda >1$ such that $|T^{'}(x) |>1$ for all $x$ where defined. 

(5)  Bounded Distortion: There exist constants $0<C_1\le C_2$ such that for all $A \in \bigvee_{j=0}^n T^{-j}\cP$
and all $x,y \in A$,
\[
C_1\le | \frac{DT^{n} (x)}{DT^{n} (y)}|\le C_2
\]

A Gibbs-Markov map  $T$ has exponential decay in $\BV(X)$, meaning that there
are $\lambda \in (0,1)$, $C>0$ such that the transfer operator
$P:L^1(\mu)\to L^1(\mu)$ defined by
\[
  \int_X f\circ T \cdot g \d \mu = \int_X f \cdot P(g) \d \mu, \quad \text{
    for all } f\in L^\infty(\mu), g\in L^1(\mu)
\]
satisfies
\begin{equation}\label{eq:P_exp-decay}
  \|P^k(g)\|_{\BV} \le C \lambda^k   \|g\|_{\BV}, \quad \text{ for } g \in \BV(X) \text{
    with } \int_X g \d \mu = 0, \text{ and } k\ge  0
\end{equation}

\section{Intermittent Maps.}\label{IM}

Here we consider a simple class of intermittent type maps
$T_{\gamma} : [0,1] \rightarrow [0,1]$, which we will call LSV maps as
defined by~\cite{LSV99}, given by
  \begin{equation}\label{IM}
    T_{\gamma} (x) := \left\{
      \begin{array}{ll}
         (2^{\gamma}x^{\gamma} +1)x& \mbox{if $0 \leq x<\frac{1}{2}$};\\
        2x-1 & \mbox{if $\frac{1}{2} \le x \le 1$}.
      \end{array} \right.
  \end{equation}
For $\gamma \in [0,1)$, there is a unique absolutely continuous ergodic invariant probability measure $\mu_{\gamma}$ with density $h_{\gamma}$ bounded away from
zero and satisfying $h_{\gamma}(x) \sim Cx^{-\gamma}$ for $x$ near zero. The existence of  stable laws, and moreover  the  existence of functional limit theorems or weak invariance principles for H\"{o}lder functions on  this class of intermittent maps has been thoroughly examined in~\cite{Gouezel_Intermittent}.
 For instance, when $\gamma\in [0,1/2)$, $n^{-1/2}\sum_{i=0}^{n-1} \phi\circ T_{\gamma}^j$ follows the  CLT ; when $\gamma = 1/2$ and $\phi(0)\neq 0$,  $(n\log n)^{-1/2}\sum_{i=0}^{n-1} \phi \circ T_{\gamma}^j$ follows the CLT; when $\gamma \in (1/2, 1)$ and $\phi(0)\neq 0$, $n^{-\gamma}\sum_{i=0}^{n-1} \phi \circ T_{\gamma}^j$ follows a stable law where  the index is $ \gamma^{-1}$.  
 Gou\"ezel~\cite[Theorem 1.3]{Gouezel_Intermittent} gives the characteristic function of the stable law for normalized  H\"older $\phi$ as 
 \[
 e^{-c |t|^{\frac{1}{\gamma}} (1-\beta \sign t \tan(\pi/2\gamma))}
 \]
 where $\beta=\sign \phi(0)$ and 
 \[
 c=\frac{h_{\gamma} (1/2)}{4 \gamma^{\gamma^{-1}}}  \phi(0)^{\gamma^{-1}}\Gamma(1-\frac{1}{\gamma})\cos (\pi/2\gamma)
 \]
 The dependence of the characteristic function on only $\phi(0)$ and $h_{\gamma} (1/2)$ is explained by the fact that the stable law for $\phi$ 
 may be  obtained by
 inducing (and then lifting) the constant function  $\phi(0)$ on the usual Young Tower for $T_{\gamma}$ with base $[1/2,1]$.

 In this paper in the setting of LSV maps we consider "bad" observables,
 for example $\phi (x)=d(x,x_0)^{-\frac{1}{\alpha}}$. Our result in this
 setting is Corollary~\ref{thm:SL-T-gamma}. For an observable $\phi$ which behaves locally as 
 $d(x,x_0)^{-\frac{1}{\alpha}}$ close to a point $x_0\not =0$ and is H\"older elsewhere  one expects a 
 competition between the stable law coming from the slow-mixing property of
 the LSV map if $\gamma \in (1/2,1)$ and the stable law arising from the
 tail of the unbounded observable $\phi$. One technical issue that arises immediately is to prove the convergence to a 
 stable law for $\phi$ in a slowly mixing system. A natural technique to try is to induce, prove that the induced system satisfies a
 stable law and then lift.  If $\frac{1}{\alpha}\ge\gamma$ this approach works in a straightforward manner. 
 Furthermore  if $\gamma>\frac{1}{\alpha}$ and $\phi (0)-E[\phi]\not =0$ 
 then a stable law of index $\frac{1}{\gamma}$ holds for a restriction of the observable in the neighborhood of the
 indifferent fixed point. This effect dominates and in fact we obtain the same stable law with index $\frac{1}{\gamma}$ we would
 obtain if $\phi$ were H\"older with $\phi (0)-E[\phi]\not =0$ i.e. with the same formula for $\beta$ and $c$ above with $\phi (0)$ replaced by
 $\phi(0)-E[\phi]$.

 However suppose $1<\alpha<2$ and   $\phi$ is locally of form $d(x,x_0)^{-\frac{1}{\alpha}}$, H\"older elsewhere, with $\phi(0)-E[\phi]=0$, 
 for example with  $E[\phi]=0$
 and  with support bounded away from the indifferent fixed point. In this setting if  $\gamma >\frac{1}{\alpha}$ 
 we are only able to prove we may lift in the 
 parameter range $\alpha < 1+\frac{1}{\gamma^2}-\frac{1}{\gamma}$.  
 If this condition holds we show that   the stable law of index $\alpha$ dominates and we obtain Birkhoff sum
  convergence to a stable law of index $\alpha$. This latter results relies on a form of the law of the iterated logarithm valid for this parameter range~\cite{DGM}.

  Finally we note that the case of $\phi(x)=d(x,0)^{-\frac{1}{\alpha}}$ has been
 clarified by Gou\"ezel, and here the two effects combine so that a stable law holds with   scaling constants  $b_n=n^{\gamma+\frac{1}{\alpha}}$
 if $1/2<\gamma+\frac{1}{\alpha}<1$.

\section{Stable limits for Birkhoff sums of dependent variables.}\label{sec:Ancilliary_Results}

Our results are based upon the investigations and results of
R.~Davis~\cite{Davis} and R.~Davis and T.~Hsing~\cite{DH95} into the partial
sum convergence of dependent random variables with infinite variance.

We paraphrase~\cite[Theorem 3.1]{DH95} below.

\begin{prop}[{\cite[Theorem 3.1]{DH95}}]\label{thm:Davis-Hsing}
  Let $ \{X_j\}$ be a 
  stationary sequence of  random variables on a probability space $(X,\mu)$ such that:
 \begin{enumerate}
 \item [(i)]
   \[
     n \mu (\frac{X_1}{b_n} \in \cdot ) \to_v \nu (\cdot)
   \]
   where
   \[
     \nu (dx)=[p\alpha x^{-\alpha-1}\one_{(0,\infty)} + (1-p)\alpha
     (-x)^{-\alpha-1}\one_{(-\infty,0)}]dx
   \]
   and $\to_v$ denotes vague convergence on $\R\setminus\{0\}$; and
 
 \item [(ii)]
   \[
     N_n:=\sum_{j=1}^n \delta_{{X_j}/{b_n}} \to_d N =
     \sum_{i=1}^{\infty}\sum_{j=1}^{\infty} \delta_{P_iQ_{ij}}
   \]
   where the convergence is in the space of random counting measures, 
   $\sum_{i=1}^{\infty}\delta_{P_i}$ is a Poisson process with intensity
   measure $\nu$, $Q_i:=\sum_{j=1}^{\infty}\delta_{Q_{ij}}$, $i\ge 1$, are
   point processes that are iid, $Q_{ij}\in [-1,1]/\{0\}$, and all point processes are mutually
   independent.
 \end{enumerate}

 Then:
 \begin{itemize}

 \item [(a)] For $0<\alpha <1$,
 \[
 \frac{1}{b_n} \sum_{j=1}^{n}X_j \to_{d} S
 \]
 where $S=\sum_{i=1}^{\infty} \sum_{j=1}^{\infty} P_i Q_{ij} $ has a stable
 distribution with index $\alpha$.
 
\item [(b)] If $1\le \alpha <2$ and
 \begin{equation}\label{eq:condition-NC}
    \qquad \lim_{\epsilon\to 0} \limsup_{n\to \infty} P \{ |
   \frac{1}{b_n} \sum_{j=1}^n \Xjtruncatedeps - \frac{1}{b_n}
   E [\sum_{j=1}^n \Xjtruncatedeps]|>\delta\}=0
   \quad \text{ for all $\delta >0$,}
 \end{equation}
 then
 \[
   \frac{1}{b_n} \left(\sum_{j=1}^{n} X_j - E [ \sum_{j=1}^{n}
     \Xjtruncatedone]\right) \to_{d} S
 \]
 where $S$ is the distributional limit of
 \[
   \sum_{i=1}^{\infty}\sum_{j=1}^{\infty} P_i Q_{ij} \one_{(\epsilon,
     \infty)} (|P_iQ_{ij}|) -\int_{\epsilon <|x| \le 1} x\nu (dx)
 \]
 as $\epsilon \to 0$. $S$ has a stable distribution with index $\alpha$.

 \end{itemize}

\end{prop}

\begin{rmk}\label{rmk:usual-stable-description}
  Condition (i) above is equivalent to
  \begin{equation}\label{eq:slow-variation1}
    {P(|X_1|> x)} =x^{-\alpha} L(x)
  \end{equation}
  and
  \begin{equation}\label{eq:slow-variation2}
    \lim_{x\to\infty} \frac{P(X_1 > x)}{P(|X_1|> x)} = p
  \end{equation}
  for a slowly varying function $L(x)$ and $0\le p \le 1$. See
  \cite[Introduction]{DH95}.
\end{rmk}

\begin{rmk}\label{RM-Che}
 
  By Chebyshev's inequality, Condition~\eqref{eq:condition-NC} is implied
  by
  \begin{equation}\label{eq:DH3.2}
    \lim_{\epsilon \to 0} \limsup_{n\to \infty} E
    (|\frac{1}{b_n}\sum_{j=0}^{n-1} \Xjtruncatedeps -
    \frac{1}{b_n} E(\sum_{j=0}^{n-1} \Xjtruncatedeps)|^2)=0, 
  \end{equation}

  By~\cite[Theorem 3]{Davis}, \eqref{eq:DH3.2} is implied 
   by
  \begin{equation}\label{eq:D9}
    \lim_{\epsilon \to 0} \limsup_{n\to
      \infty}\frac{n}{b_n^2}\sum_{j=1}^{n} \max\{ 0, E(Y_1 Y_j)\}=0,
  \end{equation}
  where $Y_j = \Xjtruncatedeps - E (\Xjtruncatedeps)$.
\end{rmk}

\begin{rmk}
  Marta-Tyran Kaminska's work~\cite[Theorem 1.3]{TK-dynamical} has the same
  condition, Equation~\eqref{eq:condition-NC},  in the case $1<\alpha \le 2$, but requires convergence in
  (ii) to a simple Poisson process i.e. $Q_{ij}=1$ for $i=j=1$ and $0$
  otherwise. Her condition was motivated by the goal of establishing
  functional limit theorems rather than distributional convergence of
  Birkhoff sums. 
\end{rmk}

\section{Main Results}\label{Main_Results}

\begin{thm}\label{thm:mix_NC}
  Suppose $(T,X,\mu)$ is a
   Gibbs-Markov map of the unit interval $X=[0,1]$.
  Let $\phi: X\to \R$ be in the domain of attraction of a stable law of
  index $\alpha \in (1,2)$ and suppose there exists a constant $K$
  such that $\|\phi \one_{|\phi|<t}\|_{BV}\le Kt$ for all large $t$.
  Define $b_n$ as in Definition~\ref{def:scaling-constants}, by
  $\lim_{n \to \infty} n \mu(|\phi | > b_n ) = 1.$ Then for all $\delta>0$,
  \begin{equation}
    \qquad \lim_{\epsilon\to 0} \limsup_{n\to \infty} \mu \{ |
    \frac{1}{b_n} \sum_{j=1}^n \phi\circ T^j \one_{\{|\phi\circ T^j |<b_n\epsilon\}}  - \frac{1}{b_n}
    E [\sum_{j=1}^n \phi\circ T^j \one_{\{|\phi\circ T^j |<b_n\epsilon}\} |>\delta\}=0
 \end{equation}

  \end{thm}
  
  \begin{rmk}
  
  The condition $\|\phi1_{|\phi|<t}\|_{BV}\le Kt$ for all large $t$ is satisfied, e.g., if
$\phi$ has finitely many intervals of monotonicity.
 For  example we are able to verify  condition~\eqref{eq:condition-on-observable} for observables such as $\phi(x) = 3 |x - x_1|^{ - 2/3}-6 |x - x_2|^{ - 2/3}$
where $x_1,x_2\in [0,1]$.

\end{rmk}


Although Theorem~\ref{thm:mix_NC} holds  for
functions $\phi:X\to \R$  that satisfy the following condition
\begin{equation}\label{eq:condition-on-observable}
  \text{There is a constant $K>0$ such that for $t$ sufficiently large
    $\bvnorm{\phi\cdot \one_{\{|\phi|<t\}}} \le K t$}
\end{equation}
we will restrict now to observables of form $\phi(x)=d(x,x_0)^{-\frac{1}{\alpha}}$ where $\alpha \in (0,2]$. 
This is because we rely on recent work~\cite{FFM20} which has shown that for such observables on  Gibbs-Markov maps the corresponding  counting process $N_n$
converges to a  Poisson point process, and this   is  key to verifying the conditions of ~\cite[Theorem 3.1]{DH95}.


\bigskip

We now state what is basically a corollary to Theorem~\ref{thm:mix_NC}.

\begin{cor}\label{thm:main_SL}
  Suppose $(T,X,\mu)$ is a
   Gibbs-Markov map of the unit interval $X=[0,1]$.
    Let $\phi(x)=d(x,x_0)^{-\frac{1}{\alpha}}$ where $\alpha \in (0,2]$, $x_0\in (0,1)$.
   Define $b_n$ as in Definition~\ref{def:scaling-constants}, by
  $\lim_{n \to \infty} n \mu(|\phi | > b_n ) = 1.$

  \begin{enumerate}
  \item [(a)] If $0<\alpha <1$ then
    \[
      \frac{1}{b_n} \sum_{j=1}^{n}\phi \circ T^j \to_{d} S
    \]
    where $S=\sum_{i=1}^{\infty} \sum_{j=1}^{\infty} P_i Q_{ij} $ has a
    stable distribution with index $\alpha$.
 
  \item [(b)] If $1\le \alpha <2$ then
      \begin{equation}\label{eq:stable-alpha-ge1}
      \frac{1}{b_n} \left(\sum_{j=1}^{n} \phi  \circ T^j - E
        [\sum_{j=1}^{n} \phi  \circ T^j \one_{\{|\phi  \circ T^j
          |<b_n\}}]\right) \to_{d} S
    \end{equation}
    where $S$ is the distributional limit of
    \[
      \sum_{i=1}^{\infty}\sum_{j=1}^{\infty} P_i Q_{ij} \one_{(\epsilon,
        \infty)} (|P_iQ_{ij}|) -\int_{\epsilon <|x| \le 1} x\nu (dx)
    \]
    as $\epsilon \to 0$. $S$ has a stable distribution with index $\alpha$.
  \end{enumerate}

  \end{cor}

\begin{rmk}\label{rmk:truncation}

In the case $0<\alpha<1$ the result is a consequence of~\cite{FFM20} who proved convergence of the counting point process $N_n$ (the conclusion of Theorem~\ref{thm:mix_NC} is not needed in this case).  Theorem~\ref{thm:mix_NC} combined with the results of~\cite{FFM20} imply the conclusion in case (b).

            
\end{rmk}

We now give an application to intermittent-type maps, describing the interplay between the slow-mixing parameter $\gamma$ and
the heavy tails parameter $\alpha$.

\begin{thm}\label{thm:SL-T-gamma}
  Suppose $(T_{\gamma} ,X,\mu)$ is a LSV map of the unit interval and
  $0\le \gamma <1$. Suppose $\phi (x) =d(x,x_0)^{-\frac{1}{\alpha}}$  where $x_0\in (0,1]$. 
  If $\alpha\in [1,2)$ and $\phi$ is integrable then we define $c_n=E[\phi]=\int d(x,x_0)^{-\frac{1}{\alpha}}d\mu_{\gamma}$, otherwise $c_n=0$.
    \begin{itemize}

  \item [(i)] If $\frac{1}{\alpha}\ge \gamma$ then

      \[
        \frac{1}{n^{1/\alpha} }\sum_{j=1}^{n}[\phi \circ T^j -c_n] \to_{d} S
      \]
      where $S=\sum_{i=1}^{\infty} \sum_{j=1}^{\infty} P_i Q_{ij} $ has a
      stable distribution of index $\alpha$;
      
      \item [(ii)] if $\frac{1}{\alpha} < \gamma $ and $\phi (0)-E[\phi] \not =0$ then
         \[
        \frac{1}{n^{\gamma}} \sum_{j=1}^{n}( \phi \circ T^j - E [ \phi])
        \to_{d} S
      \]
      where $S$ has a stable distribution of index $\gamma$. 

    \item [(iii)] if    $\phi(0)-E[\phi]=0$ and  $\frac{1}{\gamma} < \alpha <1+\frac{1}{\gamma^2}-\frac{1}{\gamma}$ then
      \[
        \frac{1}{n^{1/\alpha}}  \sum_{j=1}^{n} (\phi \circ T^j - E[\phi])
           \to_{d} S
      \]
      where $S$ has a stable distribution with index $\alpha$.
   
    \end{itemize}

  
\end{thm}
  
\begin{rmk}
 To satisfy $\gamma> \frac{1}{\alpha}$ in  case (ii) and case (iii) above it is necessary that $\alpha \in (1,2)$. The extra condition, $\frac{1}{\gamma} < \alpha <1+\frac{1}{\gamma^2}-\frac{1}{\gamma}$,  in case (iii) occurs because
 we rely on a result of Dedecker, Gou\"ezel and Merlev\`ede~\cite{DGM} which is given in the Appendix. Our standard `lifting' argument fails in this case and we rely on a 
 law of the iterated logarithm result in ~\cite{DGM} which is known to hold in this parameter regime. 

\end{rmk} 

\begin{rmk}
In~\cite[Section 2.2.1]{FFT20} it is shown that (i) holds for $\gamma\in (0,0.289)$ and $0<\alpha<1$ (actually they prove a stronger functional
convergence in the $F^{'}$ topology which implies a stable law) and it is conjectured that convergence in $F^{'}$ holds for $0<\alpha<1$  and all $\gamma <\frac{1}{2}$.
\end{rmk} 
\begin{rmk}
  The case where $\phi$ is a function of the distance to the origin $0$ has
  been clarified by Gou\"ezel~\cite{Gouezel_Intermittent}. In the set-up of
  the LSV maps where $0\le \gamma <1$ if
  $\phi (x)=x^{-\frac{1}{\alpha}}$, (so that $x_0=0$) and
  $1>\frac{1}{\alpha} +\gamma >\frac{1}{2}$ then $\phi$ converges to a
  stable law in distribution and the corresponding scaling constant is
  $n^{\gamma+\frac{1}{\alpha}}$. If $\frac{1}{\alpha} +\gamma <\frac{1}{2}$
  then we have a CLT.
\end{rmk}

\section{Proof of Theorem \ref{thm:mix_NC}.}

Recall the Karamata estimates of Proposition~\ref{prop:karamata} for
regularly varying functions.

\begin{rmk}
  Although we consider the case of a Gibbs-Markov map $T:X\to X$ and
  $\phi(x):= d(x,x_0)^{-1/\alpha}$, we are using only the following (e.g.,
  no need for the Markov property):

  \begin{itemize}
  \item for the map $T:X \to X$, $X\subset \R$:
    \begin{itemize}
    \item big images w.r.t. the invariant measure 

    \item uniform expansion: there is $\theta \in (0,1)$ such that
      $|T'(x)| \ge \theta^{-1}$ for each $x$ where the derivative exists


    \item exponential decay on $\BV$ of the transfer operator $P$ of $T$
      w.r.t. the invariant measure $\mu$ 


    \item bounded distortion

    \item invariant measure comparable to Lebesgue: density bounded above,
      and away from zero
    \end{itemize}

  \item for the observation $\phi$ that
    \begin{equation}\label{eq:condition-on-observable-bis}
      \text{There is a constant $K>0$ such that $\bvnorm{\phi\cdot
          \one_{\{|\phi|<t\}}} \le K t$ for $t$ sufficiently large.}
  \end{equation}
\end{itemize}

Condition~\eqref{eq:condition-on-observable-bis} is satisfied, e.g., if $\phi$
has finitely many intervals of monotonicity.

\end{rmk}

\begin{proof} [\textbf{Proof of Theorem \ref{thm:mix_NC}}]

  Let $([0,1],\cB,\mu,T,\cP)$ be an expanding Gibbs-Markov system as in
  Section~\ref{GM}. We will check the hypotheses of
  Theorem~\ref{thm:Davis-Hsing}.

Condition (i) is satisfied since $\phi$ is in the domain of attraction of a
stable law of index $\alpha$ (see
Remark~\ref{rmk:usual-stable-description}).

Condition (ii) holds by~\cite{FFM20}. 
Recall that by \eqref{eq:tail1} and
\eqref{eq:slow-variation1}
\[
  n\sim 1/\mu (|\phi|> b_n) = b_n^\alpha L(b_n)^{-1}
\]
Since $L$ grows slower than any power (see Lemma~\ref{thm:growth}), we will
sometimes abuse notation and consider that
\begin{equation}\label{eq:b_n-vs-n}
  b_n\sim n^{1/\alpha} 
\end{equation}


\textbf{Consider now the case of $\alpha\in(1,2)$.}

We need to establish~\eqref{eq:condition-NC}.
By Remark~\ref{RM-Che}, condition~\eqref{eq:condition-NC} is implied by
\begin{equation}\label{eq:D9bis}
  \lim_{\epsilon \to 0} \limsup_{n\to
    \infty}\frac{n}{b_n^2}\sum_{j=1}^{n} \max\left\{0, \int \Phi_n \cdot \Phi_n\circ T^{j}
    d\mu\right\}=0,
\end{equation}
where, for a fixed $\epsilon>0$, we denote
\begin{equation*} 
  \phi_n :=\phi\cdot \one_{\left\{|\phi| \le \epsilon b_n\right\}}
  \text{ and } \Phi_n := \phi_n - E (\phi_n).
\end{equation*}

To obtain~\eqref{eq:D9bis}, by the exponential decay of
correlations~\eqref{eq:P_exp-decay}, we need only show that
\begin{equation}\label{eq:k-log-n}
  \lim_{\epsilon \to 0} \limsup_{n\to
    \infty}\frac{n}{b_n^2}\sum_{j=1}^{\lfloor k \log n
    \rfloor}\max\left\{0, \int \Phi_n\Phi_n\circ T^{j}d\mu\right\}=0, 
\end{equation}
where $k$ is independent of $n$ and $\epsilon$.

Since $\mu$ is $T$-invariant, can rewrite the covariance
$\int \Phi_n\Phi_n\circ T^{j}d\mu$ as
$E(\phi_n \cdot \phi_n\circ T^{j})-[E (\phi_n)]^2$. Because
$\phi\in L^1(\mu)$, one can neglect the $[E(\phi_n)]^2$ terms in~\eqref{eq:k-log-n} as
their contributions is of order
\[
  \frac{n}{b_n^2} ( E(\phi_n))^2 \log{n} \le \frac{n}{b_n^2}
  ( E(|\phi|))^2 \log{n} \sim ( E(|\phi|))^2 n^{1-\frac{2}{\alpha}}
  \log{n}
\]
and $\alpha < 2$.

 Thus, it suffices
 to show
\begin{equation}\label{eq:D9sub}
  \lim_{\epsilon \to 0} \limsup_{n\to
    \infty}\frac{n}{b_n^2}\sum_{j=1}^{\lfloor k \log n \rfloor}\int
  |\phi_n| \cdot |\phi_n|\circ T^{j} d\mu=0.
\end{equation}

Introduce
\begin{equation*}
  \frac {1}{2} < \psi < 1, \quad u_n:=b_n^\psi,  \quad U_n := \{|\phi| \ge u_n\}.
\end{equation*}
Since $\phi$ is in the domain of attraction of a stable law with index
$\alpha$ (see \eqref{eq:slow-variation1} in
Remark~\ref{rmk:usual-stable-description}),
\begin{equation*}
  \mu (U_n) =  u_n^{-\alpha} L(u_n).
\end{equation*}

From Karamata's Theorem~\ref{prop:karamata}
and \eqref{eq:slow-variation1} we have
\begin{equation}\label{eq:computation1}
  \int {\phi_n}^{2} d\mu = \int \phi^2 \cdot \one_{|\phi|\le \epsilon b_n}
  d\mu \sim \frac{\alpha}{2-\alpha} (\epsilon b_n)^2  \mu (|\phi| > \epsilon
  b_n) = C_\alpha \epsilon^2 b_n^{2} (\epsilon b_n)^{-\alpha} L(\epsilon b_n)
\end{equation}
\begin{equation}\label{eq:computation2}
  \int_{U_n^c} {\phi}^{2} d\mu \sim \frac{\alpha}{2-\alpha} u_n^2
  \mu (|\phi| > u_n) = C_\alpha b_n^{2\psi} b_n^{-\alpha\psi} L(u_n)
\end{equation}
We decompose the sum of integrals in \eqref{eq:D9sub} as
${\rm (I)} + {\rm (II)} + {\rm (III)}$, where
\[ {\rm (I)} = \sum_{j= 1}^{\lfloor k \log n \rfloor} \int_{U_n \cap
    T^{-j}U_n} |\phi_n| \cdot |\phi_n|\circ T^{j} d\mu,
\]
\[ {\rm (II)} = \sum_{j= 1}^{\lfloor k \log n \rfloor} \int_{U_n \cap
    T^{-j}U_n^c} |\phi_n| \cdot |\phi_n|\circ T^{j} d\mu
\]
and
\[ {\rm (III)} = \sum_{j= 1}^{\lfloor k \log n \rfloor} \int_{U_n^c}
  |\phi_n| \cdot |\phi_n|\circ T^{j} d\mu.
\]

Consider (II) and (III) first.

For (III), using that $\mu$ is $T$-invariant, we have
\begin{equation}\label{eq:GM_III}
  \begin{aligned} 
    \int_{U_n^c} |\phi_n| \cdot |\phi_n| \circ T^j d\mu 
    \le \left( \int_{U_n^c} {\phi}^{2} d\mu \right)^{\frac 1 2}
    \left( \int \phi_n^2 \circ T^j d\mu \right)^{\frac 1 2} 
    = \left( \int_{U_n^c} {\phi}^{2} d\mu \right)^{\frac 1 2}
    \left(\int \phi_n^2  d\mu \right)^{\frac 1 2}.
  \end{aligned}
\end{equation}
Similarly, for (II),
\begin{equation}\label{eq:GM_II}
  \begin{aligned}
    \int_{U_n \cap T^{-j}U_n^c} |\phi_n| & \cdot
    |\phi_n| \circ T^{j} d\mu
     \le \left( \int \phi_n^2 d\mu \right)^{\frac 1 2}
      \left( \int_{T^{-j}U_n^c} {\phi}^{2} \circ T^{j} d\mu
      \right)^{\frac 1 2}\\
    & =
      \left( \int \phi_n^2 d\mu \right)^{\frac 1 2}
      \left( \int {(\phi}^{2} \cdot \one_{U_n^c})\circ T^jd\mu
      \right)^{\frac 1 2}
      =
      \left( \int \phi_n^2 d\mu \right)^{\frac 1 2}
      \left( \int {\phi}^{2} \cdot \one_{U_n^c}d \mu
      \right)^{\frac 1 2}\\
    &
      =
      \left( \int \phi_n^2 d\mu \right)^{\frac 1 2}
      \left( \int_{U_n^c} {\phi}^{2} d \mu
      \right)^{\frac 1 2}
  \end{aligned}
\end{equation}
By \eqref{eq:computation1} and \eqref{eq:computation2}
we obtain
\begin{equation*}
  \left(\int \phi_n^2 d\mu \right)^{\frac 1 2}
  \left( \int_{U_n^c} {\phi}^{2} d \mu
  \right)^{\frac 1 2} \le C_\alpha \epsilon^{1-\frac{\alpha}{2}}
  b_n^{(1-\frac{\alpha}{2})(1+\psi)} L(\epsilon b_n)^{1/2} L(b_n^\psi)^{1/2}
\end{equation*}
By \eqref{eq:tail1} and \eqref{eq:slow-variation1},
\[
  n\sim 1/\mu (|\phi|> b_n) = b_n^\alpha L(b_n)^{-1}
\]
which gives
\begin{equation*}
  \frac{n}{b_n^2}[{\rm (II)} + {\rm (III)}] \le 2C_\alpha k \epsilon^{1-\frac{\alpha}{2}}
  b_n^{-(1-\frac{\alpha}{2})(1-\psi)}\log{n} \cdot \left(\frac{L(\epsilon b_n)
      L(b_n^\psi)}{L(b_n)^2}\right)^{1/2} 
\end{equation*}
Since $L$ is slowly varying, it grows slower than any power (see
Lemma~\ref{thm:growth} in the Appendix), so, because $\psi < 1$,

\begin{equation}\label{eq:II+III}
  \limsup_{n\to \infty}\frac{n}{b_n^2} [{\rm (II)} + {\rm (III)}] =0
\end{equation}

It remains to bound (I).

Denote by $\{A^{(m)}_t\}_{t \ge 1}$ the partition induced by
$\bigvee_{j=0}^{m-1}T^{-j}\cP$.

Consider some fixed $1\le j \le k\log n$; in order to estimate
$\int_{U_n \cap T^{-j}U_n} |\phi_n| \cdot |\phi_n|\circ T^{j} d\mu$, we
have the following three possibilities.

{\bf Case 1:} $U_n\subset A^{(j)}_r$ for some $r\in \bN$.

Using the H\"older inequality and the expression of the transfer operator
$P$, 
\[
  a_j:= \int_{U_n \cap T^{-j}U_n} |\phi_n| \cdot |\phi_n|\circ T^{j} d\mu
  \le \left(\int_{X} \phi_n^2 d\mu \right)^{1/2} \left(\int_{U_n}
    \phi_n^2\circ T^{j} d\mu\right)^{1/2} = \left(\int_{X} \phi_n^2
    d\mu\right)^{1/2} \left(\int_{X} P^j(\one_{U_n}) \phi_n^2\right)^{1/2}
\]
with
\[
  P^j(\one_{U_n})\vert_x = \frac{h(y)}{h(x)} \cdot \frac {1}{(T^j)'(y)}
\]
where $y\in U_n \subset A^{(j)}_r$ is the unique point such that
$T^j(y)=x$, and $h$ is the density of the invariant measure,
$\d \mu = h \d \Lebesgue$, bounded above and away from zero.
Since $T$ is piecewise expanding, we obtain that
\[
  \|P^j(\one_{U_n})\|_{L^\infty(X)} \le C \theta^{j}
\]
for $C>0$ independent of $j$ and $n$. Thus
\begin{equation*}
  a_j\le C \theta^j  \int \phi_n^2  \d\mu
\end{equation*}

{\bf Case 2:} $U_n\subset A^{(j)}_r \cup A^{(j)}_{r+1}$ for some
$r\in \bN$.

Consider $U_n\cap A^{(j)}_r$ and $U_n\cap A^{(j)}_{r+1}$. They both satisfy
{\bf Case 1}, and therefore we have
\[
  b_j:=\int_{U_n}|\phi_n(x)| |\phi_n(T^jx)| \d\mu \le 2 C \theta^j \int
  \phi_n^2 \d\mu
\]

{\bf Case 3:} $A^{(j)}_r\subset U_n$ for some $r\in \bN$.

There exists $r_1,r_2\in \bN$ such that $A^{(j)}_{r_1}, A^{(j)}_{r_2}$
cover the endpoints of $U_n$, therefore, by {\bf Case~1},
\[
  c_j:=\int_{U_n\cap(A^{(j)}_{r_1} \cup A^{(j)}_{r_2}) }|\phi_n|
  |\phi_n|\circ T^j \d \mu \le 2 C \theta^j \int \phi_n^2 \d\mu
 \]

For the sets $A^{(j)}_r\subset U_n$, by the bounded distortion of
Gibbs-Markov system,
\[
  \mu(A^{(j)}_r\cap T^{-j}U_n) \le C \mu(A^{(j)}_r)
  \mu(U_n)/\mu(T^j{A^{(j)}_r}).
\]
Therefore, by the big image property, 
\[
  \sum_{A^{(j)}_r\subset U_n}\mu(A^{(j)}_r\cap T^{-j}U_n)\le
  \widetilde{C}\sum_{A^{(j)}_r\subset U_n}\mu(A^{(j)}_r)\mu(U_n) \le
  \widetilde{C}\mu(U_n)^2
  \]
and then
\begin{align*}
  d_j:=\sum_{\{ r : A^{(j)}_r\subset U_n\}} & \int_{A^{(j)}_r \cap (U_n \cap
                                              T^{-j}U_n)} |\phi_n||\phi_n|  \circ T^j \d \mu 
                                              \le  \sum_{\{ r : A^{(j)}_r\subset U_n\}} \int_{A^{(j)}_r\cap
                                              T^{-j}U_n}|\phi_n||\phi_n| \circ T^j \d \mu \\
  \le & \left[\sum_{\{ r : A^{(j)}_r\subset U_n\}}\mu(A^{(j)}_r\cap T^{-j}U_n)\right]
        \|\phi_n\|_{L^\infty}^2 \le  C \mu(U_n)^2 \|\phi_n\|_{L^\infty}^2
\end{align*}

We now collect all these estimates.

Using Karamata's estimate \eqref{eq:computation1} of $\bE(\phi_n^2)$, the choice of
$b_n$ given by \eqref{eq:tail1}, and the expression of
$\mu(U_n)=\mu(|\phi| > u_n)$ given by \eqref{eq:slow-variation1}:
\[
  \begin{aligned}
    \frac{n}{b_n^2} (I) \le & \frac{n}{b_n^2} \sum_{j=1}^{\lfloor k\log{n}
      \rfloor} [a_j+b_j+c_j + d_j] \le C \frac{n}{b_n^2}
    \sum_{j=1}^{\lfloor k\log{n} \rfloor} \left[\theta^j \bE(\phi_n^2) +
      \mu(U_n)^2 \|\phi_n\|_{L^\infty}^2\right]\\ & \le C
    \frac{n}{b_n^2}\left[ \epsilon^{2} b_n^{2} P(|\phi|\ge \epsilon b_n) +
      u_n^{-2\alpha}L(u_n)^2 (\epsilon b_n)^2 \log n\right]\\
    & = C[\epsilon^2 n P(|\phi|\ge \epsilon b_n) + n b_n^{-2\alpha \psi}
    \epsilon^2 L(b_n^\psi)^2 \log n] \to C \epsilon^2 \text{ as } n \to
    \infty
  \end{aligned}
\]
because $\psi > 1/2$ and $L$ grows slower than any power,
Lemma~\ref{thm:growth}.

Together with \eqref{eq:II+III}, this shows that condition \eqref{eq:D9sub}
is satisfied.
\end{proof}

\section { Proof of Theorem~\ref{thm:SL-T-gamma}}

Tyran-Kaminska~\cite[Theorem 4.4]{TK-dynamical} has
proved convergence to a simple Poisson process in our setting of Gibbs-Markov maps if $\tau (|\phi| >\epsilon b_n)\to \infty$ for 
all $\epsilon >0$, where $\tau$ is the return time function.  This non-recurrence condition is not satisfied if   $\phi$ is maximized at a periodic point. 
However recently the complete convergence of $N_n$ to a  Poisson process has  been established
in the case of $\phi(x)=d(x,x_0)^{-\frac{1}{\alpha}}$  if $x_0$ is periodic~\cite{FFM20}.  These two results cover all cases as shown 
by a dichotomy result in~\cite{FFM20}.

We will induce and model the system as a Young tower over a Gibbs-Markov base map.
As  $x_0$  need not be contained in
$[1/2,1]$ we may need to induce over a base larger than the usual Young
Tower base of $[1/2,1]$ used  for the LSV map.

Let $T_L$ denote the left branch of $T$. We consider the partition of
$(0,1]$ into  sets $(A_i)$ and $(B_j)$. We define $A_i\subset [1/2,1]$ to be that set of points in $[1/2,1)$ where  the first return time to $[1/2, 1]$ under $T$ is
$i$ and then define  $B_j=[T_{L}^{-j-1} (1/2), T_{L}^{-j} (1/2)]\subset (0,1]$,
$j\ge 0$.  Note that the sets $(A_i)$ constitute the usual partition of the base $[1/2,1)$ for the usual Young tower for the 
LSV map but we will adjoin some of the  sets $(B_j)$. Since $x_0\not =0$ there
exists a minimal $M$ such that $x_0\in  [1/2,1]\cup (\cup_{j=1}^M B_j)$.
 Define $Y:= [1/2,1]\cup (\cup_{j=1}^M B_j)$. Inducing on $Y$ the
return map to $Y$ is a Gibbs-Markov map (though not necessarily with full
branches). We take a  Tower model for $(T_{\gamma}, X,\mu)$ as
a tower over $Y$ with countable partition of the base $Y$ consisting of $(A_i)$ and
$(B_j)$ in $Y$.  If $R(x)$ is the first return to $Y$ the $T^{R} A_i=B_M$
for all $A_i \subset [0,1]$. If $2\le j \le M$ then $T^{R} B_j=B_{j-1}$ and $T^{R} B_1=[1/2,1]$.
The map $F:=T^{R}: Y\to Y$ is a Gibbs-Markov map, though not necessarily with full branches. In the 
case that $x_0 \in [1/2,1)$  we may take $F: Y\to Y$ to be a full-branched Gibbs-Markov map. We define $R_n (x)=R(x)+R(Fx)+\ldots + R(F^{n-1} x)$.

The induced map $F=T^{R}:Y\to Y$ has an invariant probability measure $\mu_Y$, whose density is Lipschitz and bounded away from infinity and $0$. 
Define $\bar{R}=\int_Y Rd\mu_Y=\frac{1}{\mu(Y)}$ by Kac's lemma.

We begin with Case (ii).

\subsection*{Case (ii):  $\gamma > \frac{1}{\alpha}$ and $\phi(0)-E[\phi]\not =0$}.

The assumption that $\gamma >\frac{1}{\alpha}$ implies that $\alpha \in (1,2)$ and hence
$E[\phi]<\infty$. Note also that the assumption that $\gamma >\frac{1}{\alpha}$ excludes the case  $0\le \gamma \le \frac{1}{2}$.

We decompose $\phi$ as $\phi=\phi_1+\phi_2$ with
$\phi_1(x) :=\phi(x)\one_{Y^c}$ and $\phi_2 (x):=\phi (x)\one_{Y}$. 
Note that $\phi_1-E[\phi_1]$  induces in a good way  on the base $Y$. 
In fact the induced version of $\phi_1-E[\phi_1]$
lies in the Banach space of functions to which the results of ~\cite[Theorem 1.2]{Gouezel_Intermittent} apply.

Following~\cite{Melbourne_Zweimuller} we will write $\phi_1-E[\phi_1]=( \phi(0) -E[\phi_1])-\frac{1}{\mu(Y)}( \phi(0) -E[\phi_1])\one_{Y}+\psi$
where $\psi$ is defined by this equation.  Note that $E[\psi]=0$,  $\psi(0)=0$ and $\psi$ is piecewise H\"older. Thus  $\psi$ satisfies a CLT and so its Birkhoff sum converges to zero
in distribution under any scaling $b_n=n^{\kappa}$, $\kappa>\frac{1}{2}$. Hence the effect of  $\psi$ is negligible, as a scaling by $n^{\gamma}$ or $n^{1/\alpha}$ will ensure that the 
scaled Birkhoff sum $\psi$ converges in distribution to zero.

 Note that $g:=( \phi(0) -E[\phi_1])-\frac{1}{\mu(Y)}( \phi(0) -E[\phi_1])\one_{Y}$ has expectation zero, $E[g]=0$. The function $g$ induces
 the function 
 \[
 \Phi_1=(\phi(0)-E[\phi_1])(R(x)-\bar{R})
 \]
 on $Y$. 
 For $x\in Y$,
 \[
 \sum_{j=0}^{\bar{R}n} g\circ T^j=\sum_{j=0}^n \Phi_1\circ F^j  +V_n (x)
 \]
 where, if $\bar{R}n\ge R_n(x)$,
 \[
 V_n (x)=\sum_{R_n(x)}^{\bar{R}n} g\circ T^j
 \] 
 and if  $\bar{R}n < R_n(x)$
 \[
 V_n (x)=-\sum_{\bar{R}n}^{R_n(x)} g\circ T^j
 \] 
 Thus we have
 \[
 \sum_{j=0}^{\bar{R}n} g\circ T^j=\sum_{j=0}^n \Phi_1\circ F^j  + \sum_{n=0}^{\bar{R}n} \psi\circ T^j +V_n (x)
 \]
 \begin{equation}\label{first}
 =(\phi(0)-E[\phi_1](R_n (x)-n\bar{R})+V_n (x) + \sum_{n=0}^{\bar{R}n} \psi\circ T^j
 \end{equation}
 We will use this observation when considering the induced form of $\phi_2-E[\phi_2]$.


We induce the observable $\phi_2$ on the Gibbs-Markov base $Y$ by defining 
$\Phi_2(x)=\sum_{i=0}^{R(x)-1}\phi_2 \circ T^{i} (x)$ where $R$ is
the first return time to $Y$ under $T$. Since $\phi_2$ has support in $Y$, $\Phi_2=0$ on all levels of 
the tower except for the base level, identified with $Y$, and on $Y$ we have $\phi_2=\Phi_2$. 

$\Phi_2$ is in the domain of 
attraction of a stable law of index $\alpha$ on the probability space $(Y, \mu_Y)$ and $E_{\mu_Y} [\phi_2]=E[\phi_2]/\mu(Y)$. 
Note that for large $t$, $\mu_Y (\Phi_2>t)=\frac{1}{\mu(Y)} \mu (\phi>t)$ and hence the  $b_n$ scaling for $\Phi_2$ is $(n\bar{R})^{\frac{1}{\alpha}}$.

 From our result on Gibbs-Markov maps $\Phi_2$
satisfies a stable law with index $\alpha$ under $F:=T^{R}$ with scaling $(n\bar{R})^{\frac{1}{\alpha}}$. 
By our main theorem, Corollary~\ref{thm:main_SL}
\[
(n\bar{R})^{-{\frac{1}{\alpha}} }\sum_{j=1}^n (\Phi_2 \circ F^j -\bar{R} E[\phi_2])\dto X_{\alpha}
\]


We write
\[
\sum_{j=0}^{[\bar{R}n]} (\phi_2 \circ T^j -E[\phi_2])= \sum_{j=0}^{R_n (x) } (\phi_2 \circ T^j -E[\phi_2]) +W_n (x)
\]
where $W_n(x)=\sum_{R_n (x)}^{[\bar{R}n]}(\phi_2 \circ T^j-E[\phi_2]) $ or $W_n(x)=-\sum^{R_n (x)}_{[\bar{R}n]}(\phi_2 \circ T^j-E[\phi_2])$.

Furthermore 
\[
\sum_{j=0}^{R_n (x) }( \phi_2 \circ T^j -E[\phi_2])=\sum_{j=0}^n(\Phi_2\circ F^j-\bar{R}E[\phi_2])
\]
\[
-E[\phi_2][R_n(x)-n\bar{R}]
\]
Thus 
\begin{eqnarray}\label{second}
\sum_{j=0}^{[\bar{R}n]} (\phi_2 \circ T^j -E[\phi_2])&=& \sum_{j=0}^n(\Phi_2\circ F^j-\bar{R} E[\phi_2])\\\nonumber
&-&E[\phi_2][R_n(x)-n\bar{R}] +W_n(x)
\end{eqnarray}

Adding Equations~\eqref{first} and \eqref{second} we obtain the representation
\begin{eqnarray}\label{decomposition}
\sum_{j=0}^{[\bar{R}n]} (\phi \circ T^j -E[\phi])&=& \sum_{j=0}^n(\Phi_2\circ F^j-\bar{R} E[\phi_2]) +\sum_{j=0}^{\bar{R}n} \psi\circ T^j \\
&+&(\phi(0)-E[\phi_1]-E[\phi_2])[R_n(x)-n\bar{R}] +V_n (x) +W_n(x)\nonumber
\end{eqnarray}

As noted before $n^{-\kappa} \sum_{j=0}^{\bar{R}n} \psi\circ T^j $ converges in distribution to zero for any $\kappa>\frac{1}{2}$.

We will show that 
\[
\frac{1}{(\bar{R}n)^{\gamma}} V_n(x)\dto 0
\]
and 
\[
\frac{1}{(\bar{R}n)^{\gamma}} W_n(x)\dto 0
\]
in distribution. This implies that 
\[
(\bar{R}n)^{-\frac{1}{\alpha}} \sum_{j=0}^{[\bar{R}n]} (\phi \circ T^j -E[\phi])\dto (\bar{R}n)^{-\frac{1}{\alpha}}\sum_{j=0}^n(\Phi_2\circ F^j-\bar{R} E[\phi_2])
\]
and hence 
\[
(\bar{R}n)^{-\frac{1}{\alpha}} \sum_{j=0}^{[\bar{R}n]} (\phi \circ T^j -E[\phi])\dto X_{\alpha}
\]
\

Now we show 
\[
\frac{1}{(\bar{R}n)^{\gamma}} W_n(x)\dto 0
\]
in distribution. The proof for $V_n (x)$ is the same mutatis mutandis.

Since  $\gamma >\frac{1}{2}$
\[
\frac{R_n-n\bar{R}}{n^{\gamma}}\dto X_{\frac{1}{\gamma}}
\]
as the return time function $R$ lies in the domain of attraction of $X_{\frac{1}{\gamma}}$ and satisfies the conditions of~\cite[Theorem 1.2]{Gouezel_Intermittent}.
Given $\epsilon >0$ choose $L$ large enough that for all $x$ in a set $G_0$, $\mu(G_0^c)<\epsilon$
\[
 |\frac{R_n (x)-n\bar{R}}{n^{\frac{1}{\gamma}}}|<L
\]

As a consequence of the ergodic theorem, given $\epsilon>0$, there exists an $M_1$ such that for a set $G_1$, $\mu(G_1^c)<\epsilon$,
\[
|\frac{1}{m}\sum_{j=0}^{Lm} (\phi_2\circ T^j (x) -E[\phi_2])|<\epsilon
\]
for all $m\ge M_1$ and all $x\in G_1$. Note that this implies that for all $M_1\le k \le Lm$,
\[
|\frac{1}{m}\sum_{j=0}^{k} (\phi_2\circ T^j (x) -E[\phi_2])|<\epsilon
\]
As $\phi_2$ is integrable we may choose $M_2>M_1$ such that for all $x\in G_2\subset G_1$ with
$\mu(G_1 \setminus G_2)<\epsilon$
\[
|\frac{1}{M_2}\sum_{j=0}^{M_1} (\phi_2\circ T^j (x) -E[\phi_2])|<\epsilon
\] 
Hence for all $x\in G_2$, for all $m>M_2$, for all $0\le k \le Lm$

\[
|\frac{1}{m}\sum_{j=0}^k (\phi_2\circ T^j (x) -E[\phi_2])|<\epsilon
\]
Since also for all $m>M_2$
\[
|\frac{1}{m}\sum_{j=0}^{Lm} (\phi_2\circ T^j (x) -E[\phi_2])|<\epsilon
\]
we have for all $0\le k \le Lm$
\[
|\frac{1}{m}\sum_{j=k}^{Lm} (\phi_2\circ T^j (x) -E[\phi_2])|<2\epsilon
\]
Now choose $n$ large enough that $n^{\gamma}>M_2$. Note that $\mu(T^{-\bar{R}n} (G_0\cap G_2))=\mu(G_0 \cap G_2)$ and $\mu(T^{-(\bar{R}n-n^{\gamma}L)} (G_0\cap G_2))=\mu (G_0\cap G_2)$.

If $T^{\bar{R}n} x\in G_2\cap G_0$ and $R_n(x)>\bar{R}n$ then  $R_n (x)-\bar{R}n< n^{\gamma}L$
and hence $n^{-\gamma} |\sum_{\bar{R}n}^{R_n(x)} \phi_2\circ T^j -E[\phi_2]|<2\epsilon$.
Similarly  if $T^{\bar{R}n-n^{\gamma}L} x\in G_2\cap G_0$ and $R_n(x)<\bar{R}n$ then  
$\bar{R}n -R_n (x)< n^{\gamma}L$ and hence $n^{-\gamma} |\sum_{R_n(x)}^{n\bar{R}} \phi_2\circ T^j -E[\phi_2]|<2\epsilon$. This shows that $\frac{1}{n^{\gamma}}W_n(x)
\dto 0$ in distribution. This completes the proof of case (ii).

\subsection*{Case (i) : The case $\frac{1}{\alpha}\ge \gamma$.} 

We suppose $1<\alpha<2$  and define $E_{Y}[\Phi_2]$ where the expectation is on $(Y,\mu_Y)$, so that $E_Y[\Phi_2]=\bar{R} E[\phi_2]$.
The argument we give works equally well for $0<\alpha<1$ by taking $E_Y[\Phi_2]=0$. We will leave out the dependence on $Y$ and 
write simply $E[\Phi_2]$.

From our result on Gibbs-Markov maps 
\[
(\bar{R}n)^{-\frac{1}{\alpha}} \sum_{j=0}^{n}[\Phi_2\circ F^j-E[\Phi_2]] 
\]
converges to a stable law $X_{\alpha}$ of index $\alpha$. We will use
~\cite[Theorem 4.6]{Gouezel_Doubling} (see Appendix~\ref{sec.gouezel-lifting}) to
lift this stable law to a stable law for $\phi_2$ under $T$. We verify
condition (b) of Proposition~\ref{thm.gouezel-lifting} of the Appendix to
show that
\[
n^{-\frac{1}{\alpha}} \sum_{j=0}^{n}(\phi_2 \circ T^j -\mu (\phi_2))
\]
converges in distribution to $X_{\alpha}$. In Condition (b) of Proposition~\ref{thm.gouezel-lifting} of the Appendix we take $\alpha(n)=n^{\gamma}$, $C_n=nE[\Phi_2]$ and $B_n=(n\bar{R})^{\frac{1}{\alpha}}$.

$R$ is integrable on the probability space $(G, \mu_G)$ with expectation $\bar{R}$. We define
$R_n (x)=R(x)+R(Fx)+\ldots+R(F^{n-1} x)$. By the ergodic theorem
\[
\lim_{n\to \infty} \frac{R_n (x)-n\bar{R}}{n}=0
\]
for $\mu_G$ a.e. $x\in G$. Now we show that $R$ satisfies a stable law of index $\frac{1}{\gamma}$ under $F$ (this result
is well-known).

The return-time function $R$ is constant on partition elements of $G$ and hence measurable with respect to the partition on $G$. $R$ is in the domain of attraction of a stable law of index $\frac{1}{\gamma}$ if $\frac{1}{2}<\gamma <1$ or the 
central limit theorem if $\gamma<\frac{1}{2}$.  By~\cite[Corollary 4.3]{TK-dynamical} 

\[
\frac{1}{n^{\gamma}} \sum_{j=0}^{n} [R\circ F^j-n\bar{R}]
\]
 converges to a stable law of index $\frac{1}{\gamma}$.

Hence
\[
n^{-\gamma} \sum_{j=0}^{n} [R\circ F^j-n\bar{R}]
\]
is tight.

By~\cite[Theorem 4.6]{Gouezel_Doubling} 
this implies $\phi_2$ satisfies a stable law of index $\alpha$ and scaling
$b_n=n^{\frac{1}{\alpha}}$ and centering $E(\phi_2)/\bar{R}$ (which is $E(\Phi_2)$).

Now   $\phi_1$ satisfies either a CLT (if $\phi_1(0)-E[\phi]=0$) or a stable law with scaling $n^{\gamma}$ (if $\phi_1 (0)-E[\phi_1] \not =0$), and thus  the scaled  Birkhoff sum
$n^{-\frac{1}{\alpha}} \sum_{j=0}^{n}\phi_1 \circ T^j -E[\phi_1] $ converges in distribution to zero. This proves that 
\[n^{-\frac{1}{\alpha} } \left(\sum_{j=1}^{n} (\phi_1 \circ T^j - E(\phi_1)) + (\phi_2 \circ T^j - E(\phi_2)) \right)
\]
\[
        =n^{-\frac{1}{\alpha} } \left(\sum_{j=1}^{n} (\phi \circ T^j - E(\phi) \right) \dto_{d} X_{\alpha}
      \]
          where $X_{\alpha}$ has a stable distribution with index $\alpha$ given by~\cite{DH95}.

\subsection*{Case (iii):  $\gamma>\frac{1}{\alpha} $, $\phi(0)-E[\phi]=0$ and $ \alpha<1+\frac{1}{\gamma^2}-\frac{1}{\gamma}$} 

Recall the representation~\eqref{decomposition}, 
\begin{eqnarray*}
\sum_{j=0}^{[\bar{R}n]} (\phi \circ T^j -E[\phi])&=& \sum_{j=0}^n(\Phi_2\circ F^j-\bar{R} E[\phi_2]) +\sum_{j=0}^{\bar{R}n} \psi\circ T^j \\
&+&(\phi(0)-E[\phi_1]-E[\phi_2])[R_n(x)-n\bar{R}] +V_n (x) +W_n(x)
\end{eqnarray*}
As before $n^{-\frac{1}{\alpha}}\sum_{j=0}^{\bar{R}n} \psi\circ T^j $ converges to zero in distribution  and under our assumption $\phi(0)-E[\phi]=0$, we have 
\begin{eqnarray*}
(\bar{R}n)^{-\frac{1}{\alpha}}\sum_{j=0}^{[\bar{R}n]} (\phi \circ T^j -E[\phi])&\dto & (\bar{R}n)^{-\frac{1}{\alpha}}\sum_{j=0}^n(\Phi_2\circ F^j-\bar{R} E[\phi_2]) 
 + (\bar{R}n)^{-\frac{1}{\alpha}} [V_n (x) +W_n(x)]
\end{eqnarray*}

We will use a result of Dedecker, Gou\"ezel and Merlev\`ede~\cite[Theorem 1.7]{DGM} (see Appendix) to obtain the almost sure convergence of
\[
n^{-\frac{1}{\alpha}}\sum_{j=0}^{n^{\gamma}}(\phi_2\circ T^j -E[\phi_2])
\]
and 
\[
n^{-\frac{1}{\alpha}}\sum_{j=0}^{n^{\gamma}}(\phi_1\circ T^j -E[\phi_1])
\]
to zero  in the parameter range we consider. In fact the proof in the case of $\phi_1-E[\phi_1]$ is the same as that 
of  $\phi_2-E[\phi_2]$, so we give only the latter proof.

 This
is the key ingredient in the proof of this section that allows us to control the discrepancy $W_n(x)$ and 
to lift the induced stable law of index $\alpha$ from $\Phi_2\circ F^j $ to $\phi_2\circ T^j$ under the corresponding scalings.

There is a law of the iterated
logarithm for Birkhoff sums satisfying an exact (not asymptotic) stable
law~\cite{Chover} (for example iid random variables in the domain of attraction of a stable law)
which unfortunately is not applicable in our setting. 
In 
\[
n^{-\frac{1}{\alpha}}\sum_{j=0}^{n^{\gamma}}(\phi_2\circ T^j -E_{\mu}[\phi_2] )
\]
let $m=n^{\gamma}$ (leaving out integer part notation) then we may rewrite as 
\[
m^{-\frac{1}{\alpha\gamma}}\sum_{j=0}^{m}(\phi_2\circ T^j -E_{\mu}[\phi_2] )
\]

For any sufficiently small $\delta>0$  we will show that we  may take $p=\alpha\gamma+\delta$ in 
Dedecker, Gou\"ezel and Merlev\`ede~\cite[Theorem 1.7]{DGM} and conclude
\[
\lim_{n\to \infty} n^{-\frac{1}{\alpha}}\sum_{j=0}^{n^{\gamma}}(\phi_2\circ T^j -E_{\mu}[\phi_2]) =0
\]
$\mu$ a.e.

We need to check the conditions of~\cite[Theorem 1.7]{DGM}.

Our function $\phi_2$ falls in their class  of functions
$\mathcal{F}(H,\mu)$ where the tail function is $H(t) \sim t^{-\alpha}$ as $t\to \infty$.  For small $\delta>0$
the condition
$1< p \le 2$ and $\gamma < \frac{1}{p}$  are satisfied if $\gamma^2<\frac{1}{\alpha}$ as $p=\alpha\gamma+\delta$.
Now we consider condition (1.7)
\[
H(t)^{(1-p\gamma)/(1-\gamma)}\le C t^{-p}
\]
This condition is satisfied if $p<\frac{\alpha}{\alpha\gamma+1-\gamma}$. Taking $\delta>0$ small this condition follows
if $\gamma<\frac{1}{\alpha\gamma+1-\gamma}$, which is equivalent to 
 $\gamma-\gamma^2<1-\alpha \gamma^2$. The condition  $\gamma-\gamma^2<1-\alpha \gamma^2$ 
imposes more restrictions than $\gamma^2<\frac{1}{\alpha}$. Thus the conditions
of Dedecker, Gou\"ezel and Merlev\`ede~\cite[Theorem 1.7]{DGM} are satisfied
in our setting if
$\frac{1}{\gamma}< \alpha<1+\frac{1}{\gamma^2}-\frac{1}{\gamma}$. As an illustrative example, if $\gamma=\frac{2}{3}$
we require $\frac{3}{2}< \alpha <\frac{7}{4}$.

By the ergodic theorem
\[
\lim_{n\to \infty} \frac{1}{n}\sum_{j=0}^{n} [\phi_2\circ T^j (x) -E_{\mu}[\phi_2] ] =0~(\mu~a.e. x \in X)
\]
Since $[n\bar{R}-R_n (.)]/n^{\gamma}$ converges in distribution to a stable law of index $\frac{1}{\gamma}$ given $\epsilon>0$
we may choose $L$  and $M_1$ large enough that
\[
\mu_G \{ x \in G: | n\bar{R}-R_n (x) |>Tn^{\gamma} \}<\epsilon
\] 
for all $n\ge M_1$.

Thus for all $n\ge M_1$ the set
\[
B_n:=\{ x: |n\bar{R}-R_n (x)|>\tau n^{\gamma}
\]
satisfies $\mu (B_n)<\epsilon$.

In our parameter range
\begin{eqnarray}\label{LIL}
  \lim_{n\to \infty} n^{-\frac{1}{\alpha}} \sum_{j=0}^{n^{\gamma}}( \phi_2 \circ T^j -E[\phi_2]) =0
\end{eqnarray}
for $\mu$ a.e. $x\in [0,1]$. 

Hence 
\begin{eqnarray}\label{LIL}
  \lim_{n\to \infty} n^{-\frac{1}{\alpha}} \sum_{j=0}^{n^{\gamma}L}( \phi_2 \circ T^j -E[\phi_1]) =0
\end{eqnarray}
for $\mu$ a.e. $x\in [0,1]$.

Choose $M_2>M_1$ large enough that 
\[
\mu \{ x\in X: \max_{M_2\le k \le L n^{\gamma}} |n^{-\frac{1}{\alpha}}
|\sum_{j=0}^{T(k^{\gamma})}[ \phi_2 \circ T^j (x) -E[\phi_2]]| >\epsilon \}<\epsilon
\]

Note that this implies that for all $n>M_2$,
\[
\mu \{ x\in X: \forall M_2 \le k \le L n^{\gamma}, |n^{-\frac{1}{\alpha}} \sum_{j=0}^{T(k-1)}[ \phi_2 \circ T^j (x) ]| <\epsilon \}>1-\epsilon
\]

By measure preservation 
\[
\mu_G \{ x\in G: \forall M_2 \le k \le L  n^{\gamma}, |n^{-\frac{1}{\alpha}}\sum_{j=0}^{T(k-1)}[ \phi_2 \circ T^j (T^{n\bar{R}} x)]| <\epsilon_1 \}>1-\epsilon
\]
and
\[
\mu_G \{ x\in G: \forall M_2 \le k \le L n^{\gamma}, |n^{-\frac{1}{\alpha}}\sum_{j=n\bar{R} -n^{\gamma}}^{T(k-1)}[ \phi_2 \circ T^j (T^{n\bar{R}-n^{\gamma}} x)]| <\epsilon \}>1-\epsilon_1
\]

Thus except for a set of points $x\in G$ of  $\mu_G$ measure less than $2\epsilon$
\[
| n^{-\frac{1}{\alpha}} [\sum_{j=0}^{n\bar{R}} \phi_2 \circ T^j (x) -E[\phi_2]-(n\bar{R})^{-\frac{1}{\alpha}} \sum_{j=0}^{n} \Phi_2 \circ F^j (x) -E[\Phi_2]|<2\epsilon
\]
This implies 
\[
n^{-\frac{1}{\alpha}} [\sum_{j=0}^{n\bar{R}} \phi_2 \circ T^j (x)-E[\phi_2]]
\]
converges in distribution to a stable law of index $\alpha$.

Once the
conclusion of Equation~\eqref{LIL} is established the proof proceeds as in
the previous section.  In any case the stable law
for $\Phi_2$ lifts to $\phi_2$ and if $\phi_1 (0)=0$ then the scaling is
$n^{\frac{1}{\alpha}}$, if $\phi_1 (0)\not =0$ the scaling is $n^{\gamma}$.

\qed

\section{Discussion}

In Theorem 8.1 we show that small jumps are ``negligible'' for a wide class of heavy-tailed functions on Gibbs-Markov maps. This
result is used to investigate the interplay between the effects of heavy-tails and slow-mixing in a common model of intermittency
for observables of form $\phi(x)=d(x,x_0)^{-\frac{1}{\alpha}}$. Our results in this direction rely on work of~\cite{FFM20} who
proved complete convergence of a corresponding counting point process to a Poisson process.
Our results are for stable laws but suggest that convergence in stronger topologies may hold for all $\alpha>\frac{1}{\gamma}$, $0<\gamma<1$.

\section{Appendix}

\begin{lem}\label{thm:growth}
  A slowly varying function $L$ grows slower than any power.\end{lem}
\begin{proof} Let $\delta > 0$ be arbitrary. Using the Representation
  Theorem (see e.g.~\cite[Theorem 1.3.1]{Bingham-Goldie-Teugels-1987}):
  \[
    \frac{L(x)}{x^\delta} \sim \frac{c(x)\exp\left(\int_{1}^x
        \frac{\epsilon(s)}{s} d s\right)}{\exp(\delta\int_1^x \frac{1}s \d
      s)} = {c(x)\exp\left(\int_{1}^x \frac{\epsilon(s) - \delta}{s} d
        s\right)}
  \]
  with $c(x) \to c \in (0, \infty)$ and $\epsilon(x) \to 0$ as $x\to \infty$.
\end{proof}


\subsection{A result of Gou\"ezel}\label{sec.gouezel-lifting}

We use the following result of Gou\"{e}zel~\cite[Theorem 4.6]{Gouezel_Doubling}:

\begin{prop}\label{thm.gouezel-lifting}
Let $(T,X,\mu)$  be an ergodic probability preserving map, let $\alpha (n)$ and  and $B_n$ be two sequences of integers which are regularly varying with positive indexes.  Let $A_n\in \R$ and let $Y \subset X$ be a subset with positive measure. We will denote by $\mu_Y (.):= \frac{\mu |_{Y}  }{\mu (Y)}$ the induced probability measure.
Let $R: Y \to \bN$ be the return time of $T$  to $Y$ and $F =T^{R} : Y \to Y$ be the  the induced map. Define $\bar{R}=\int_Y R d\mu= \frac{1}{\mu (Y)}$. Consider a measurable function $\phi: X\to \R$ and define $\Phi : Y \to \R$ by $\Phi (y)=\sum_{j=0}^{R(y)-1} \phi\circ T^j$.
Define $S_n (\Phi)=\sum_{j=0}^{n-1} \Phi\circ F^j$. Assume that 
\[
\frac{ S_n (\Phi)-A_n}{B_n}
\]
converges in distribution (with respect to $\mu_Y$) to a random variable $S$.

  Additionally assume that either:
  \begin{enumerate}
  \item [(a)] $\frac{\sum_{j=0}^n R\circ F^j- n\bar{R}}{\alpha(n)}$ tends
    in probability to zero and
    $\max_{0\le k\le \alpha (n)} \frac{|S_k (\Phi)|}{B_n}$ is tight
  \item []
  \item [or]
  \item []
  \item [(b)] $\frac{\sum_{j=0}^n R\circ F^j- n\bar{R}}{\alpha(n)}$ is tight and
    $\max_{0\le k\le \alpha (n)} \frac{|S_k (\Phi)|}{B_n}$ tends in probability
    to zero.
  \end{enumerate}

  Then 
  \[
    \big(\sum_{j=0}^{n-1} \phi\circ T^j -A_{\left \lfloor n\mu(Y)\right
      \rfloor}\big)/B_{\left \lfloor n\mu (Y)\right \rfloor }
  \]
  converges in distribution (with respect to $\mu$) to $S$.

\end{prop}

\subsection{A result of Dedecker, Gou\"ezel and Merlev\`ede}

We paraphrase the results of Dedecker, Gou\"ezel and Merlev\`ede~\cite[Theorem 1.7]{DGM} that we use for the benefit of the reader. They define a class
of functions  $\mathcal{F}(H,\mu)$. Let $\mu$ be a probability measure on $\R$
and $H$ a tail function. Let $\mbox{Mon}(H,\mu)$ denote the set of functions $f: \R \to \R$ which are monotonic on some open interval and null elsewhere
such that $\mu (|f|>t)\le H(t)$. We define  $\mathcal{F}(H,\mu)$ to  be the closure in $L^1(\mu)$ of the set of functions that can be written as
$\sum_{j=0}^l a_j f_j$ where $\sum_{j=0}^l |a_j | \le 1$ and $f_j \in \mbox{Mon}(H,\mu)$. In our setting if $\phi_2$ is integrable 
then $\phi_2\in   \mathcal{F}(H,\mu)$. Suppose the LSV map has parameter $0<\gamma <1$, $\phi_2$ satisfies $\mu (|\phi_2| >t)\sim t^{-\alpha}$ (and hence
$H(t)\sim t^{-\alpha}$. As a consequence of Dedecker, Gou\"ezel and Merlev\`ede~\cite[Theorem 1.7]{DGM} if:
\begin{itemize}
\item[(i)] $1<p\le 2$ and $0<\gamma <\frac{1}{p}$ 
\item[(ii)] $H(t)^{(1-p\gamma)/(1-\gamma)} \le C t^{-p}$
\end{itemize}
then for any $b>\frac{1}{p}$ 
\[
  n^{-\frac{1}{p}} (\ln (n))^{-b} \sum_{j=0}^{n-1} [\phi_2\circ T^j -\mu
  (\phi_2)] \to 0 \quad \text{ $\mu$-a.e.}
\]
In our setting we take $p=\alpha \gamma$.

\bibliographystyle{alpha} 

\bibliography{references_all}

\end{document}